\documentclass[11pt]{article}
\usepackage{latexsym, amsfonts}
\newcommand{\be}{\begin{equation}}
\newcommand{\ee}{\end{equation}}
\newcommand{\bea}{\begin{eqnarray}}
\newcommand{\eea}{\end{eqnarray}}
\newcommand{\bean}{\begin{eqnarray*}}
\newcommand{\eean}{\end{eqnarray*}}
\newcommand{\brray}{\begin{array}}
\newcommand{\erray}{\end{array}}
\newcommand{\ben}{\begin{equation}{nonumber}}
\newcommand{\een}{\end{equation}{nonumber}}

\newtheorem{dfn}{Definition}[section]
\newtheorem{thm}[dfn]{Theorem}
\newtheorem{lmma}[dfn]{Lemma}
\newtheorem{ppsn}[dfn]{Proposition}
\newtheorem{crlre}[dfn]{Corollary}
\newtheorem{xmpl}[dfn]{Example}
\newtheorem{rmrk}[dfn]{Remark}
\newcommand{\bdfn}{\begin{dfn}}
\newcommand{\bthm}{\begin{thm}}
\newcommand{\blmma}{\begin{lmma}}
\newcommand{\bppsn}{\begin{ppsn}}
\newcommand{\bcrlre}{\begin{crlre}}
\newcommand{\bxmpl}{\begin{xmpl}}
\newcommand{\brmrk}{\begin{rmrk}}
\newcommand{\edfn}{\end{dfn}}
\newcommand{\ethm}{\end{thm}}
\newcommand{\elmma}{\end{lmma}}
\newcommand{\eppsn}{\end{ppsn}}
\newcommand{\ecrlre}{\end{crlre}}
\newcommand{\exmpl}{\end{xmpl}}
\newcommand{\ermrk}{\end{rmrk}}

\newcommand{\IC}{\mathbb{C}}

\newcommand{\IN}{{I\! \! N}}

\newcommand{\IR}{\mathbb{R}}

\newcommand{\IZ}{\mathbb{Z}}

\newcommand{\al}{\alpha}

\newcommand{\cla}{{\cal A}}
\newcommand{\clb}{{\cal B}}
\newcommand{\clc}{{\cal C}}

\newcommand{\clh}{{\cal H}}

\newcommand{\clk}{{\cal K}}

\newcommand{\cln}{{\cal N}}

\def\a*{{\cal A}_{h,*}}
\def\B{{\cal B}(h)}
\def\B1{{\cal B}_1(h)}
\def\b{{\cal B}^{s. a. }(h)}
\def\b1{{\cal B}^{s. a. }_1(h)}

\newcommand{\ot}{\otimes}

\newcommand{\raro}{\rightarrow}

\def \qed {$\Box$}

\begin{document}
\begin{center}
{\large {\bf Invariance and the twisted  Chern
character : a case study  }}\\
{\bf Debashish Goswami}\\
{\bf Stat-Math Unit, Indian Statistical Institute}\\
{203, B. T. Road, Kolkata 700108, India.}\\
E-mail : goswamid@isical.ac.in \\
\end{center}

\begin{abstract}
We give details of the proof of the remark made in \cite{G2} that the Chern
characters of the canonical generators on the K homology of the quantum group
$SU_q(2)$ are not invariant under the natural $SU_q(2)$ coaction.
Furthermore, the conjecture made in \cite{G2} about the nontriviality of the twisted
Chern character coming from an odd equivariant spectral triple on $SU_q(2)$ is
settled in  the affirmative.
\end{abstract}

\section{Introduction}
Noncommutative geometry (NCG) (a la Connes, see \cite{Con} ) and the $C^*$-algebraic
 theory of
quantum groups (see, for example, \cite{Wo1}, \cite{Wo2})  are two well-developed
mathematical areas which share the basic idea of `noncommutative
mathematics', namely, to view a general (noncommutative) $C^*$
algebra as noncommutative analogue of a topological space, equipped
with additional structures resembling
 and generalizing those in the classical (commutative) situation,
 e.g. manifold or Lie group structure. A lot of fruitful interaction
 between these two areas is thus quite expected. However, such an
 interaction was not very common until recently, when a systematic
 effort by a number of mathematicians  for understanding  $C^*$-algebraic quantum groups
   as noncommutative manifolds in the sense of Connes triggered a
   rapid and interesting development to this direction. However,
   quite surprisingly,
   such an effort
   was met with a number of obstacles
   even in the case of the simplest non-classical quantum group, namely
   $SU_q(2)$ and it was not so clear for some time whether this
   (and other standard examples of quantum groups) could be nicely
   fitted into the framework of Connes' NCG (see \cite{G} and the
   discussion and references therein). The problem of finding a nontrivial equivariant spectral triple for
    $SU_q(2)$ was finally settled in
   the affirmative in the papers by Chakraborty and Pal (\cite{CP}, see also \cite{connew} and \cite{landietal}
    for subsequent development), which increased the hope for
    a happy marriage between NCG and quantum group theory. However,
    even in the case of $SU_q(2)$, a few puzzling questions remain
    to be answered. One of them is the issue of invariance of the
    Chern character, which we have addressed in \cite{G2} and
    attempted to suggest a solution through the twisted version of
    the entire cyclic cohomology theory, building on the ideas of
    \cite{KMT}. In that paper, we also made an attempt to study the
    connection between twisted and the conventional NCG following a
    comment in \cite{connew}. The present article is a follow-up of
    \cite{G2}, and we mainly concentrate on $SU_q(2)$, considering it as a test-case for comparing the twisted
     and conventional formulation of NCG.

\section{Notation and background}
Let $\cla=SU_q(2)$ (with $0<q<1$)  denote the
$C^*$-algebra generated by two elements $\al, \beta$ satisfying
$$\alpha^*\alpha+\beta^*\beta=I,\quad\alpha\alpha^*+q^2\beta\beta^*=I,\quad \alpha\beta-q\beta\alpha=0,$$
 $$\alpha\beta^*-q\beta^*\alpha=0,\quad\beta^*\beta=\beta\beta^*.$$
 We also denote the $\ast$-algebra generated by $\al$ and $\beta$
 (without taking the norm completion) by $\cla^\infty$. There is a
 Hopf $\ast$ algebra structure on $\cla^\infty$, as can be seen
 from, for example, \cite{Wo2}. We denote the canonical coproduct on
 $\cla^\infty$ by $\Delta$.
We shall also use the so-called Sweedler convention, which
 we briefly explain now. For $a \in \cla^\infty,$ there are
 finitely many elements $a^{(1)}_i, a^{(2)}_i$, $i=1,2,...,p$ (say), such that
 $\Delta(a)=\sum_i a^{(1)}_i \ot a^{(2)}_i.$ For notational
 convenience, we abbreviate this as $\Delta(a)=a^{(1)} \ot
 a^{(2)}.$  For any positive integer $m,$ let
 $\cla^\infty_m$ be the $m$-fold algebraic tensor product of $\cla^\infty.$
 There is a natural coaction of $\cla^\infty$ on $\cla^\infty_m$
 given by
 $$ \Delta^m_\cla (a_1 \ot a_2 \ot ...\ot a_m):=({a_1}^{(1)} \otimes ...
{a_m}^{(1)}) \otimes ({a_1}^{(2)}...{a_m}^{(2)}),$$
 using the Sweedler notation, with summation being implied. Let us
 recall the convolution $\ast$ defined in \cite{G2}. If $\phi :
 \cla^\infty_m \raro \IC$ is an $m$-linear functional, and $\psi :
 \cla^\infty \raro \IC$ is a linear functional, we define their
 convolution $\phi \ast \psi : \cla^\infty_m \raro \IC$ by the following :
 $$(\phi \ast \psi)(a_1\ot...\ot a_m):= \phi(a_1^{(1)} \ot ... \ot
 a_m^{(1)})\psi(a_1^{(2)}...a_m^{(2)}),$$ using the Sweedler
 convention. We say that an $m$-linear functional $\phi$ is {\it invariant}
  if $\phi \ast \psi=\psi(1) \phi $ for every functional $\psi$ on $\cla^\infty$.

 In   \cite{MNW}, the K-homology $K^*(\cla^\infty)$ has
been
 explicitly  computed. It has been shown there that
 $K^0(\cla^\infty)=K^1(\cla^\infty)=\IZ$, and the Chern characters (in cyclic cohomology ) of the
  generators of
  these K-homology groups, denoted by $[\tau_{\rm ev}]$ and $[\tau_{\rm odd}]$ respectively,
    are also explicitly written down.  

\section{Main results}
\subsection{Chern characters are not invariant}
In this subsection, we  give detailed arguments for a remark made in
\cite{G2} about the impossibility of having an invariant Chern
character for $\cla^\infty$ under the conventional (non-twisted)
framework of NCG. To make the notion of invariance precise, we give
the following definition (motivated by a comment by G. Landi, which
is gratefully acknowledged).
  \bdfn
  \label{inv}
  We say that a class $[\phi ] \in HC^n(\cla^\infty)$ is invariant if there
  is an invariant $n+1$-linear functional
  $\phi^\prime$ such that $\phi^\prime$ is a cyclic cocycle and $\phi^\prime \sim \phi$ (i.e.
  $[\phi]=[\phi^\prime]$).
  \edfn
 It is easy to see that the Chern chracter $[\tau_{\rm ev}]$ cannot
 be invariant. Had it been so, it would follow from the uniqueness of the 
  Haar state (say $h$)  on $SU_q(2)$ that $\tau_{\rm ev}$ must be a scalar  multiple of $h$.   Since $\tau_{\rm
 ev}$ is a nonzero trace, it would imply that $h$ is a trace too.
 But it is known (see \cite{Wo2}) that $h$ is not a trace. 
 
 However,
 proving that $[\tau_{\rm odd}]$ is not invariant requires little
 bit of detailed arguments. We begin with the following observation.
  \blmma
  \label{lem1}
  If $\tau$ is a trace on $\cla^\infty$, i.e. $\tau \in HC^0(\cla^\infty)$, then
 we have that $$(\partial \xi) \ast \tau=\partial(\xi \ast
 \tau)$$ for every functional $\xi$ on $\cla^\infty$, where the Hochschild coboundary operator $\partial$ is defined by $$ (\partial \xi)(a,b)=\xi(ab)-\xi(ba).$$
 \elmma
 {\it Proof :}\\
 We shall use the Swedler notation. We have that for $a_0,a_1 \in
 \cla^\infty$, \bean \lefteqn{(\partial \xi \ast \tau)(a_0,a_1)}\\
 &=& (\partial \xi)(a_0^{(0)} \ot a_1^{(0)}) \tau(a_0^{(1)}
 a_1^{(1)})\\
 &=&
 \xi(a_0^{(0)}a_1^{(0)})\tau(a_0^{(1)}a_1^{(1)})-\xi(a_1^{(0)}a_0^{(0)})\tau(a_0^{(1)}a_1^{(1)})\\
 &=& \xi(a_0^{(0)}a_1^{(0)}\tau(a_0^{(1)}a_1^{(1)})
 -\xi(a_1^{(0)}a_0^{(0)})\tau(a_1^{(1)}a_0^{(1)})~~({\rm
 since}~~\tau~~{\rm is~a~trace})\\
 &=& (\xi \ast \tau)(a_0a_1)-(\xi \ast \tau)(a_1a_0)\\
 &=& \partial (\xi \ast \tau)(a_0,a_1).\eean
 \qed

 The above lemma allows us to define the multiplication $\ast$ at the level of cohomology classes. More precisely,  for 
   $[\phi] \in HC^1(\cla^\infty)$ and $[\eta] \in HC^0(\cla^\infty)$, we set $[\phi] \ast [\eta]:=[\phi \ast \eta] \in HC^1(\cla^\infty)$, which is well-defined by the Lemma \ref{lem1}. Similarly $[\eta] \ast [\phi]$ and $[\eta] \ast [\eta^\prime] $ (where $[\eta^\prime ] \in HC^0(\cla^\infty)$) can be defined. We now recall from \cite{MNW} that 
   $$ [\tau_{\rm ev}] \ast
 [\tau_{\rm ev}]=[\tau_{\rm ev}],~~[\tau_{\rm ev}] \ast [\tau_{\rm
 odd}]=[\tau_{\rm odd}] \ast [\tau_{\rm ev}]=0.$$ We also note that $\tau_{\rm ev}(1)=1$ and that
 $\tau_{\rm ev}$ is a trace, i.e. $\tau_{\rm ev}(ab)=\tau_{\rm
 ev}(ba)$.

Using this observation,  we are now in a position to  prove that the Chern character of the
generator of $K^1(\cla^\infty)$ is not an invariant class. \bthm
$[\tau_{\rm odd}]$ is not invariant.
\ethm {\it Proof :}\\
Suppose that there is $\phi \sim \tau_{\rm odd}$ such that $\phi$ is
invariant. 
Then we have $$ [\phi \ast \tau_{\rm ev} ] =[\phi] \ast [\tau_{\rm ev}]=[\tau_{\rm odd} ] \ast [\tau_{\rm ev}]=0.$$ However, since we have   
 $\phi \ast \tau_{\rm ev}=\tau_{\rm
ev}(1)\phi=\phi$ by the invariance of $\phi$, it follows that $[\phi]=[\phi \ast \tau_{\rm ev}]=0$, that is, 
 $[\tau_{\rm odd}]=0$, which is a contradiction. \qed

\subsection{Nontrivial pairing with the twisted Chern character}
As already mentioned in the introduction, in \cite{G2} we have made
an attempt to recover the desirable property of invariance by making
a departure from the conventional NCG and using the twisted
 entire cyclic cohomology. We briefly recall here some of the basic concepts from that paper and refer the reader to \cite{G2} and the references
therein for more details of this approach. We shall use the results derived in that paper wihout always giving a specific reference. 

Let us  give  the definition of twisted entire cyclic cohomology  for Banach algebras  for   simplicity, but
 note that  the theory   extends  to locally convex algebras,
which we actually need. The extension to the locally convex
algebra case follows exactly as remarked in [1, page 370]. So, let
$\cla$ be a unital Banach algebra, with $\| . \|_*$ denoting its
Banach norm, and let $\sigma$ be a continuous automorphism of
$\cla$, $\sigma(1)=1$. For $n \geq 0$, let $C^n$ be the space of
continuous $n+1$-linear functionals $\phi$ on $\cla$ which are
$\sigma$-invariant, i.e.
$\phi(\sigma(a_0),...,\sigma(a_n))=\phi(a_0,...,a_n) \forall
a_0,...,a_n \in \cla$; and $C^n =\{0\}$ for $n<0$. We define
linear maps $T_n, N_n : C^n \raro C^n$,  $ U_n : C^n \raro
C^{n-1}$ and $V_n : C^n \raro C^{n+1}$ by, $$
(T_nf)(a_0,...,a_n)=(-1)^n f(\sigma(a_n), a_0,...,a_{n-1}),
N_n=\sum_{j=0}^n T_n^j,$$ $$ (U_nf)(a_0,...,a_{n-1})=(-1)^n
f(a_0,...,a_{n-1}, 1),$$
$$(V_nf)(a_0,...,a_{n+1})=(-1)^{n+1}
f(\sigma(a_{n+1})a_0,a_1,...,a_n).$$
 Let $B_n=N_{n-1}U_n (T_n-I),$ $b_n=\sum_{j=0}^{n+1} T_{n+1}^{-j-1}V_nT_n^j$.
  Let $B,b$ be maps on the complex $C\equiv (C^n)_n$ given by $B|_{C^n}=B_n,
b|_{C^n}=b_n$. It is easy to verify (similar to what is done for
the untwisted case , e.g. in \cite{Con}) that $B^2=0$, $b^2=0$ and
$Bb=-bB$, so that we get a bicomplex $(C^{n,m}\equiv C^{n-m})$
with differentials $d_1,d_2$ given by $d_1=(n-m+1)b : C^{n,m}
\raro C^{n+1,m},$ $d_2=\frac{B}{n-m} : C^{n,m} \raro C^{n,m+1}.$
Furthermore, let $C^{e}=\{ (\phi_{2n}){n \in \IN}; \phi_{2n} \in
C^{2n} \forall n \in \IN \},$ and $C^o=\{ ( \phi_{2n+1}){n \in
\IN}; \phi_{2n+1} \in C^{2n+1} \forall n \in \IN \}$. We say that
an element $\phi=(\phi_{2n})$ of $C^e$ is a $\sigma$-twisted even
entire cochain if the radius of convergence of the complex power
series $\sum \| \phi_{2n} \|\frac{z^n}{n!}$ is infinity, where $\|
\phi_{2n} \| := \sup_{ \| a_j \|_* \leq 1} | \phi_{2n}
(a_0,....,a_{2n})|.$ Similarly we define $\sigma$-twisted odd
entire cochains, and let $C^e_\epsilon(\cla, \sigma)  $
($C^o_\epsilon(\cla,\sigma)$ respectively) denote the set of
$\sigma$-twisted even (respectively odd) entire cochains. Let
$\tilde{\partial} =d_1+d_2$ , and we have the short complex
$C^e_\epsilon(\cla, \sigma) \frac{{\stackrel{\tilde{\partial}}
{\longleftarrow}}}{{\stackrel{\longrightarrow}{\tilde{\partial}}}}
C^o_\epsilon(\cla, \sigma)$. We call the cohomology of this
complex the $\sigma$-twisted entire cyclic cohomology of $\cla$
and denote it by $H^*_\epsilon(\cla,\sigma)$.
 Let $\cla_\sigma=\{ a \in \cla : \sigma(a)=a
\}$ be the fixed point subalgebra for the automorphism $\sigma$.
There is a canonical pairing $<.,.>_{\sigma,\epsilon} :
K_*(\cla_\sigma) \times H^*_\epsilon(\cla, \sigma) \raro \IC$.
 We shall need the pairing for the odd case, which we write down :
 $$ <[u],[\psi]> \equiv < [u],[\psi]>_{\sigma,\epsilon}=\frac{1}{\sqrt{2 \pi i}} \sum_{n=0}^\infty (-1)^n \frac{n!}{(2n+1)!} \psi_{2n+1}(u^{-1},u,...,u^{-1},u),$$ where $[u] \in K_1(\cla_\sigma)$ and $[\psi] \in H^1_\epsilon(\cla,\sigma).$

\bdfn
\label{spectral}
 Let $\clh$ be a separable Hilbert space,
$\cla^\infty$ be a $\ast$ subalgebra (not necessarily complete)
of $\clb(\clh)$, $R$ be a positive (possibly unbounded) operator
on $\clh$, $D$ be a
self-adjoint operator in $\clh$ with compact resolvents such that the following hold :\\
(i) $[D,a] \in \clb(\clh)$ $ \forall a \in \cla^\infty$,\\ (ii)
$R$ commutes with $D$,\\ (iii) For any real number $s$ and $a \in
\cla^\infty$, $\sigma_s(a):=R^{-s}aR^s$  is bounded and belongs to
$\cla^\infty$. Furthermore,
 for any positive integer $n$, $\sup_{s \in [-n,n]} \| \sigma_s(a) \| <
\infty.$\\ Then we call the quadruple $(\cla^\infty, \clh, D, R)$
an odd $R$-twisted spectral data. 
   We say that
the odd twisted spectral data is $\Theta$-summable
if $Re^{-tD^2}$ is trace-class for all $t >0$.
\edfn

 Let us now recall the construction of twisted Chern character from a given odd twisted spectral data  
$(\cla^\infty, \clh, D, R)$. 
  Let  $\clb$ denote the set of all $A \in
\clb(\clh)$ for which $\sigma_s(A):=R^{-s}AR^s \in \clb(\clh) $
for all real number $s$, $[D,A] \in \clb(\clh)$ and $s \mapsto \|
\sigma_s(A) \|$ is bounded over compact subsets of the real line.
In particular, $\cla^\infty \subseteq \clb.$  We define for $n \in
\IN$ an $n+1$-linear functional $F_n$ on $\clb$ by the
formula
$$ F_n(A_0,...,A_n)=\int_{\Sigma_n}Tr( A_0 e^{-t_0 D^2}
A_1e^{-t_1D^2}...A_ne^{-t_n D^2}R)dt_0...dt_n,$$  where $\Sigma_n
=\{ (t_0,...,t_n) : t_i \geq 0, \sum_{i=0}^n t_i=1 \}$. 
 
Let us now equip $\cla^\infty$ with the locally convex topology
given by the family of Banach norms $\|.\|_{*,n},n=1,2,...$, where
$\|a\|_{*,n}:= \sup_{s \in [-n,n]} ( \| \sigma_s(a) \|+\|
[D,\sigma_s(a)] \|)$. Let $\cla$ denote the completion of
$\cla^\infty$ under this topology, and thus $\cla$ is Frechet
space. We can   construct the (twisted) Chern character in
$H^o_\epsilon(\cla,\sigma)$, where $\sigma=\sigma_1$, which
extends on the whole of $\cla$ by continuity.

 \bthm
 \label{chern}
Let  $\phi^o\equiv
(\phi_{2n+1})_n$ be defined by  $$ \phi_{2n+1}(a_0,...,a_{2n+1})=\sqrt{2i}
F_{2n+1}(
a_0,[D,a_1],...,[D,a_{2n+1}]), a_i \in \cla.$$ Then
 we have  $(b+B)\phi^o=0$,  hence  $\psi^o
\equiv ((2n+1)! \phi_{2n+1})_n \in H^o_\epsilon(\cla,\sigma).$
\ethm

We shall also need  some results from the theory of semifinite spectral triples and the corresponding JLO cocycles and index formula, as discussed in, for example, \cite{carey}.
An odd  semifinite spectral triple  is given by $(\clc, \cln, \clk, D)$, where $\clk$ is a separable Hilbert space, $\cln\subseteq \clb(\clk)$ is a von Neumann algebra with a faithful semifinite normal trace (say $\tau$),
 $D$ is a self-adjoint operator affiliated to $\cln$, $\clc$ is a $\ast$-subalgebra of $\clb(\clk)$ such that $[D, c] \in \clb(\clk)$ for all $c \in \clc$.  In the terminology of \cite{carey}, $(\cln,D)$ is also called an odd, unbounded Breuer-Fredholm module for the norm-closure of $\clc$.  It is called $\Theta$-summable if $\tau(e^{-tD^2})<\infty$ for all $t >0$. For a $\Theta$-summable semifinite spectral triple, there is a canonical construction of JLO cocycle and index theorem (see \cite{carey}), which are very similar to their counterparts in the conventional framework of NCG.

Let us now settle in the affirmative  conjecture made in \cite{G2} about the nontriviality of   the twisted Chern character of a natural twisted spectral data  obtained from the equivariant spectral triple of \cite{CP}. For reader's convenience, we briefly recall the construction of this equivariant spectral triple.
  Let us  index the space of irreducible (co-)representations of
$SU_q(2)$
  by half-integers, i.e. $n=0,\frac{1}{2},1,...;$ and index the orthonormal basis
   of the corresponding $(2n+1)^2$ dimensional subspace of $L^2(SU_q(2),h)$
    by $i,j=-n,...,n, $ instead of $1,2,...,(2n+1).$ Thus, let us consider the orthonormal
     basis $e^n_{i,j},n=0,\frac{1}{2},...; i,j=-n,-n+1,...,n$ in the notation of \cite{CP}.
        We consider any of the equivariant spectral triples constructed by
the authors of \cite{CP} and in the associated Hilbert space
$\clh=L^2(SU_q(2),h)$  define the following positive unbounded operator $R$ :
$$R(e^n_{i,j})=q^{-2i-2j}e^n_{i,j},$$
 $n=0, \frac{1}{2},,1,...;$ $i,j=-n,-n+1,...,n.$  Let us choose a
spectral triple given by the Dirac operator $D$ on $\clh$, defined
by $$ D(e^n_{i,j})= d(n,i) e^n_{i,j},$$ where $d(n,i)$ are as in
(3.12) of \cite{CP}, i.e. $d(n,i)=2n+1$ if $n=i,$ $d(n,i)=-(2n+1)$
otherwise. It can  easily be seen that $(\cla^\infty, \clh, D, R )$
is an odd $R$-twisted spectral data  and
furthermore,   the fixed point subalgebra $SU_q(2)_\sigma$ for
$\sigma(.)=R^{-1}\cdot R$ is the unital $\ast$-algebra generated by
$\beta$,  so it contains
$u=I_{1}(\beta^*\beta)(\beta-I)+I$ which can be chosen to be a generator
of $K_1(SU_q(2)) =\IZ$ (see  \cite{CP}).  It is easily seen  that the map
from $K_1(C^*(u))$ to $K_1(SU_q(2))$, induced by the inclusion
map, is an isomorphism of the $K_1$-groups (where $C^*(u)$ denotes
the unital $C^*$-algebra generated by $u$).
Thus, we can  consider the
pairing of the twisted Chern character with $K_1(C^*(u))$, and in
turn with $K_1(SU_q(2))$ using the isomorphism noted before. The
important question raised in  \cite{G2} is  whether we recover the nontrivial pairing
obtained in \cite{CP} in our twisted framework, and in what follows,
 we shall give an affirmative answer to this question. 
\bthm
The pairing between $K_1(SU_q(2)_\sigma) \cong K_1(SU_q(2))$ and the (twisted) Chern character of the above twisted spectral data coincides with the pairing between $K_1(SU_q(2))$ and the Chern character of the (non-twisted) spectral triple $(\cla^\infty, \clh, D)$. In particular, this pairing is nontrivial. 
\ethm
 {\it Proof :}\\
 
  Let $\cln$ be the von Neumann algebra in $\clb(\clh)$ generated by $\beta$ and $f(D)$ for all bounded measurable functions $f : \IR \raro \IC.$  Since $R$ commutes with both $\beta$ and $D$, it is easy to see that the functional $\cln \ni X \mapsto \tau(X):=Tr(X R)$ defines a faithful, normal, semifinite trace on the von Neumann algebra $\cln$. Moreover, $(\cln, D)$ is an unbounded $\Theta$-summable Breuer-Fredholm module for the norm-closure of the  unital $\ast$-algebra (say $\clc$) generated by $\beta$.

 Moreover, it follows from  the  fact that $R$ commutes with $D$ and $u$ that  the pairing of $[u]$ with the twisted Chern character (say $\psi^o\equiv (\psi_{2n+1})$)  coming
 from the twisted spectral data $(\cla^\infty, \clh, D, R)$ is given by  
 \bean
 \lefteqn{ <[u],[\psi^o]>}\\
 &=& \frac{1}{\sqrt{2 \pi i}} \sum_{n=0}^\infty (-1)^{n} \frac{n!}{(2n+1)!} \psi_{2n+1}(u^{-1},u,...,u^{-1},u)\\
& & \frac{1}{\sqrt{\pi}}\sum_{n=0}^\infty (-1)^{n} n!~ 
 \int_{\Sigma_{2n+1}}Tr(u^{-1} 
e^{-t_0D^2}[D,u]e^{t_1D^2}...[D, u]e^{t_{2n+1}D^2} R)dt_0...dt_{2n+1},\\
&=& \frac{1}{\sqrt{\pi}}\sum_{n=0}^\infty (-1)^{n} n!~ 
 \int_{\Sigma_{2n+1}}\tau(u^{-1} 
e^{-t_0D^2}[D,u]e^{t_1D^2}...[D, u]e^{t_{2n+1}D^2} )dt_0...dt_{2n+1}
\eean
   which is nothing but  the
  pairing between $[u] \in K_1(\clc)$ and the Breuer-Fredholm module $(\cln,D)$ mentioned before. By  Theorem 10.8 of \cite{carey} and a straightforward but somewhat lengthy calculation along the lines of index computation in \cite{CP}, we can show that  the value of this pairing is equal to  $-ind_\tau(A)\equiv -(\tau(P_A)-\tau(Q_A))$ for the following operator $A : \clh_0 \raro \clh_0$, where $\clh_0$ is the closed subspace spanned by $\{ e^n_{n,j}, n=0, \frac{1}{2}, ..., j=-n, -n+1,...,n \}$, $P_A$, $Q_A$ are the orthogonal projections onto the kernel of $A$ and the kernel of $A^*$ respectively and where $r$ is a positive integer such that $q^{2r} < \frac{1}{2} < q^{2r-2}$ :
  $$ Ae^n_{n,j}=-q^{(n+j)(2r+1)}(1-q^{2(n-j)})^r(1-q^{2(n-j-1)})^{\frac{1}{2}}e^{n+\frac{1}{2}}_{n+\frac{1}{2},j-\frac{1}{2}}+(1-q^{2r(n+j)}(1-q^{2(n-j)})^r)e^n_{n,j}.$$
  It can be verified by computations as in  \cite{CP} that $Ker(A)=\{ 0 \}$ and $Ker(A^*)$ is the one dimensional subspace spanned by the vector $\xi=\sum_{n=\frac{1}{2},\frac{3}{2},...} p_n e^n_{n,-n}$, where $p_{\frac{1}{2}}=1$ and for $n \geq \frac{3}{2},$ $$ p_n=\frac{1-(1-q^{4n-2})^r}{(1-q^{4n})^{\frac{1}{2}}(1-q^{4n-2})^r}
  ...\frac{1-(1-q^2)^r}{(1-q^4)^{\frac{1}{2}}(1-q^2)^r}.$$ Clearly, since $Re^n_{-n,n}=e^n_{n,-n},$ we have $R \xi=\xi$ and thus $$-ind_\tau(A)= \frac{1}{\| \xi \|^2}\tau( |\xi><\xi|)=\frac{1}{\| \xi \|^2}Tr(R |\xi><\xi|)=1,$$  which is the same as  the value of the pairing between $[u] \in K_1(SU_q(2))$ and the conventional Chern character corresponding to the spectral triple constructed in \cite{CP}. 
  \qed

\vspace{2mm} {\noindent} Thus we see that both the conventional
and twisted frameworks of NCG give essentially the same results
for the example we considered, namely
 $SU_q(2)$. The aparent weakness of the twisted NCG arising from the fact that the twisted cyclic
 cohomology can be paired naturally with only the K theory of the invariant subalgebra
 and not of the whole algebra, does not seem to
 pose any essential difficulty for studying the noncommutative geometric aspects of
 $SU_q(2)$, since by a suitable choice of the twisting
 operator $R$ as we did one could make sure that the K theory of the corresponding invariant subalgebra is isomorphic with the K theory of
 the whole, and also the pairing between the Chern character and the generator of the
 K theory in the twisted framework is equal to the similar pairing in the ordinary (non-twisted)
 framework of NCG. It will be  important and interesting to investigate
  whether a similar fact remains true for a larger class of quantum groups, and we hope to pursue this in the future.

\end{document}